\documentclass[11pt]{article}
\usepackage[kerning, tracking, spacing]{microtype}
\microtypecontext{spacing=nonfrench}
\usepackage[margin=1in]{geometry}
\usepackage{amsmath, amsthm}
\usepackage{amssymb}
\usepackage{mathrsfs}
\usepackage[shortlabels]{enumitem}
\usepackage{graphicx}
\usepackage{mathdots}
\usepackage{color}
\usepackage{hyperref,url}
\usepackage{tikz, tikz-cd}
\usetikzlibrary{matrix}
\usetikzlibrary{shapes}
\usetikzlibrary{arrows,decorations.markings, tikzmark}
\usepackage{mathtools}
\usepackage{ stmaryrd }
\usepackage[normalem]{ulem}
\usepackage[utf8]{inputenc}
\usepackage{xcolor}
\hypersetup{
    colorlinks,
    linkcolor={red!50!black},
    citecolor={blue!50!black},
    urlcolor={blue!80!black}
}

\pgfarrowsdeclare{bad to}{bad to}
{
  \pgfarrowsleftextend{-2\pgflinewidth}
  \pgfarrowsrightextend{\pgflinewidth}
}
{
  \pgfsetlinewidth{0.8\pgflinewidth}
  \pgfsetdash{}{0pt}
  \pgfsetroundcap
  \pgfsetroundjoin
  \pgfpathmoveto{\pgfpoint{-3\pgflinewidth}{4\pgflinewidth}}
  \pgfpathcurveto
  {\pgfpoint{-2.75\pgflinewidth}{2.5\pgflinewidth}}
  {\pgfpoint{0pt}{0.25\pgflinewidth}}
  {\pgfpoint{0.75\pgflinewidth}{0pt}}
  \pgfpathcurveto
  {\pgfpoint{0pt}{-0.25\pgflinewidth}}
  {\pgfpoint{-2.75\pgflinewidth}{-2.5\pgflinewidth}}
  {\pgfpoint{-3\pgflinewidth}{-4\pgflinewidth}}
  \pgfusepathqstroke
}

\DeclareMathAlphabet{\mymathbb}{U}{BOONDOX-ds}{m}{n}

\def\Mbar{\overline{\M}}
\def\M{\mathcal{M}}
\def\Mct{\M^{\mathrm{ct}}}
\def\Mrt{\M^{\mathrm{rt}}}

\def\CH{\mathsf{CH}}

\def\R{\mathsf{R}}

\def\QQ{\mathbb{Q}}

\usepackage{amsthm}
\theoremstyle{definition}
\newtheorem{definition}{Definition}[section]
\newtheorem{theorem}[definition]{Theorem}

\newtheorem{question}[definition]{Question}

\newtheorem*{conjecture*}{Conjecture}
\newtheorem*{question*}{Question}
\newtheorem*{theorem*}{Theorem}
\newtheorem*{faberconj}{Gorenstein Conjecture}
\newtheorem*{fzconj}{Faber-Zagier Conjecture}

\title{The tautological ring of $\M_{25}$ is not Gorenstein}
\author{Aaron Pixton}
\begin{document}

\maketitle

\begin{abstract}
  We show that the tautological ring of $\M_{25}$ is not Gorenstein, disproving a conjecture of Faber from 1993. This also implies that the Faber-Zagier relations are the complete list of tautological relations on $\M_{25}$.
\end{abstract}

\section{Introduction}

Let $g \ge 2$, and let $\M_g$ be the moduli space of smooth curves of genus $g$. Mumford \cite{Mumford} introduced the {\it tautological ring} of $\M_g$, the subring of the Chow ring $\R^*(\M_g)\subseteq \CH^*(\M_g)$ generated by the kappa classes $\kappa_i\in \CH^i(\M_g)$. (In this paper, all Chow rings are taken with rational coefficients.) Determining the structure of the tautological ring is thus equivalent to determining which polynomials in the kappa classes vanish, i.e. determining all {\it tautological relations}.

In 1993 Faber proposed a beautiful conjecture giving a complete description of the tautological ring:
\begin{faberconj}[{\cite[Conjecture 1(a)]{Faber}}]
  The tautological ring $\R^*(\M_g)$ is Gorenstein with one-dimensional socle in degree $g-2$. In other words, $\R^{>(g-2)}(\M_g) = 0, \R^{g-2}(\M_g)\cong\QQ$, and the multiplication map
  \[
    \R^d(\M_g)\times \R^{g-2-d}(\M_g)\to \R^{g-2}(\M_g)\cong \QQ
  \]
  is a perfect pairing of $\QQ$-vector spaces for all $0 \le d \le g-2$.
\end{faberconj}
Properly speaking, this is not a complete description without also including a formula for the proportionalities between all kappa polynomials of degree $g-2$, which Faber also provided \cite[Conjecture 1(c)]{Faber}. Equipped with these proportionalities and assuming the conjecture, we can effectively determine whether any given kappa polynomial is zero or nonzero by computing its pairing with every kappa monomial of complementary degree.

The parts of this conjectural description that deal with kappa polynomials of high degree were all proved within a decade. Looijenga \cite{Looijenga} proved that $\R^{>(g-2)}(\M_g) = 0$ and that $\dim_\QQ \R^{g-2}(\M_g) \le 1$, and Faber \cite{Faber-nonvanishing} proved that $\dim_\QQ \R^{g-2}(\M_g) > 0$. Getzler and Pandharipande \cite{Getzler-Pandharipande} proved the proportionalities conditional on the degree-zero Virasoro conjecture for $\mathbb{P}^2$, which was then proved by Givental \cite{Givental}. This left only the perfect pairing itself as conjectural.

Faber \cite{Faber} verified the Gorenstein conjecture for all $g\le 15$ (later extended to $g\le 23$) by introducing a geometric method to prove many tautological relations. For each $g\le 23$, he was able to obtain every relation predicted by the Gorenstein conjecture in this way, but there were gaps for $g\ge 24$. For example, with $g=24$ there was exactly one ``missing relation'', appearing in codimension $12$. There are $77$ kappa monomials in $\R^{12}(\M_{24})$ and the rank of the socle pairing is $36$, so the Gorenstein conjecture predicts a $41$-dimensional space of tautological relations there. Faber's method (and every known source of tautological relations that has been studied since then) only produces relations in $\R^{12}(\M_{24})$ spanning a $40$-dimensional space. So $\dim_\QQ \R^{12}(\M_{24}) = 36\text{ or }37$, and the Gorenstein conjecture for $g=24$ is equivalent to the vanishing of a specific kappa polynomial $P_{24}$ of degree $12$ (unique only modulo the $40$-dimensional space of known relations and up to scaling).

The case of $g=25$ is similar: there is a single ``missing relation'' in $\R^{12}(\M_{25})$ whose vanishing is equivalent to the Gorenstein conjecture. We prove that in fact this kappa polynomial $P_{25}$ is nonvanishing, and thus the Gorenstein conjecture is false.

\begin{theorem}\label{thm:main}
The Gorenstein conjecture is false for $g=25$.
\end{theorem}

The general approach taken is fairly standard. We pair $P_{25}$ with a ``detector'' class $Q_{25}\in \R^{60}(\Mbar_{25})$, where $\Mbar_{25}$ is the standard Deligne-Mumford compactification of $\M_{25}$ by stable curves. This class has the following two properties:
\begin{enumerate}[(a)]
\item $Q_{25}$ vanishes on the entire boundary of $\Mbar_{25}$, in the sense that $j^*Q_{25}=0$ for $j$ equal to any of the basic gluing maps $\Mbar_{24,2}\to\Mbar_{25}$ or $\Mbar_{h,1}\times\Mbar_{25-h,1}\to\Mbar_{25}$.
  \item If $\overline{P}_{25}\in\CH^{12}(\Mbar_{25})$ is any extension of $P_{25}$ to the boundary, then the intersection number $\int_{\Mbar_{25}}\overline{P}_{25}Q_{25}$ is nonzero.
  \end{enumerate}
  Property (a) implies that the choice of extension $\overline{P}_{25}$ in property (b) does not matter. As a consequence, these two properties together imply that $P_{25}$ is nonvanishing, since if it were a relation then we could take $\overline{P}_{25} = 0$ and obtain a contradiction. (In fact this intersection number argument also works in cohomology by a standard mixed Hodge theory argument as discussed in \cite[Appendix, Proposition 1]{PZP}, so this also disproves the Gorenstein conjecture in cohomology, not just in Chow.)

This class $Q_{25}$ that detects the nonvanishing of $P_{25}$ is obtained by a fairly involved computer search, using the Faber-Zagier/$3$-spin relations \cite{PPZ,Janda} to obtain the boundary vanishing of property (a). The space $\R^{60}(\Mbar_{25})$ is extremely large, so a simplification is needed to make this search feasible, even by computer. That simplification is that we only search for classes that are multiples of $\lambda_g\lambda_{g-2}$, where $\lambda_i = c_i(\mathbb{E})$ are the Chern classes of the Hodge bundle on $\Mbar_g$. (In comparison, any multiple of $\lambda_g\lambda_{g-1}$ automatically satisfies boundary vanishing but will not detect anything beyond what the regular socle pairing does.) The details of this search (and thus the proof of Theorem~\ref{thm:main}) are given in Section~\ref{sec:computations}.

In Section~\ref{sec:background} we discuss more background and history relating to Faber's Gorenstein conjecture, including the analogous Gorenstein conjectures for other moduli spaces of curves and the Faber-Zagier relations. In Section~\ref{sec:computations} we discuss the computations involved in the construction of $Q_{25}$. Finally, in Section~\ref{sec:questions} we list a few questions for future work.

\subsection*{AI usage}
The main idea of this paper (searching for a multiple of $\lambda_g\lambda_{g-2}$ that vanishes on the boundary) was proposed and implemented by GPT-6 Astra after being prompted by the author (in Codex, ``max'' effort setting) over a period of several days of suggesting various perspectives on the problem. Probably the author's most important contribution in the prompts was to suggest searching for Gorenstein counterexamples in $\M_{20,1}$ or $\M_{25}$, not just $\M_{24}$ (where the approach does not seem to work, as discussed in Section~\ref{sec:computations}). The author then checked the arguments (including independently writing code to check the computations from scratch) and wrote this paper himself, without use of AI for any writing.

The GPT-6 Astra usage fit within the weekly usage cap of a \$100/month Pro subscription.

\subsection*{Acknowledgments}
I am grateful to Carel Faber and Rahul Pandharipande for first teaching me about the Gorenstein conjecture and the Faber-Zagier relations when I was a graduate student, and for many discussions about them in the years since.

This work was supported by the NSF grant DMS-2301506.

\section{Background}\label{sec:background}
In this section we briefly describe some of the prior evidence against the Gorenstein conjecture for $\M_g$, some of which will also serve as useful background for the computations in the following section.

\subsection{Hodge integrals and Gorenstein conjectures}
Let $\mathbb{E}$ be the Hodge bundle on $\Mbar_g$, the rank $g$ vector bundle given by pushing forward the relative dualizing sheaf of the universal curve over $\Mbar_g$. Let $\lambda_i = c_i(\mathbb{E})$ be the Chern classes of $\mathbb{E}$. Then $\lambda_g\lambda_{g-1}$ vanishes on the boundary of $\Mbar_g$; this follows from the fact that the total Chern class $c(\mathbb{E})$ pulls back to any boundary stratum of $\Mbar_g$ to the product of a copy of itself on each component, along with Mumford's identity $c(\mathbb{E})c(\mathbb{E}^*) = 1$ \cite{Mumford}.

As a consequence, the socle evaluation isomorphism
\[
  \R^{g-2}(\M_g)\to\QQ
\]
can be constructed by integrating an extension to the boundary against $\lambda_g\lambda_{g-1}$. The proportionalities thus can be interpreted as evaluating Hodge integrals
\[
  \int_{\Mbar_g}\lambda_g\lambda_{g-1}K,
\]
where $K$ is a kappa polynomial.

Faber and Pandharipande \cite{Faber-Pandharipande} noted that $\lambda_g$ itself vanishes on the boundary divisor of stable curves with a non-separating node and thus can similarly be integrated against to define a map
\[
  \CH^{2g-3}(\Mct_g)\to\QQ,
\]
where $\Mct_g$ is the complement of this boundary divisor, i.e. the moduli space of compact type curves. Motivated by this, they speculated that the tautological subring $\R^*(\Mct_g)\subseteq \CH^*(\Mct_g)$ (which they defined) might also be Gorenstein, with socle in degree $2g-3$. Together with the case of $\Mbar_g$ itself (where no Hodge classes are needed to integrate against for a socle in degree $3g-3$), this led to a triad of Gorenstein conjectures, for the tautological rings of $\M_g,\Mct_g,\Mbar_g$. In fact one can also consider all of these conjectures with marked points (using the rational tails moduli space $\Mrt_{g,n}$ in place of smooth curves) and ask whether the tautological rings of $\Mrt_{g,n},\Mct_{g,n},\Mbar_{g,n}$ are Gorenstein rings with socles of degrees $g-2+n, 2g-3+n, 3g-3+n$ respectively.

But Petersen and Tommasi \cite{Petersen-Tommasi,Petersen} proved that the Gorenstein conjecture is false for $\Mbar_{2,20}$, and Petersen \cite{Petersen} proved that it is false for $\Mct_{2,8}$. More recently, Canning \cite{Canning} proved that it is false for $\Mbar_{g,n}$ whenever $g\ge 2$ and $2g+n \ge 24$ (so in particular it is false for $\Mbar_{12}$). And Canning, Larson, and Schmitt \cite{Canning-Larson-Schmitt} proved that it is false for $\Mct_{g,n}$ whenever $g\ge 2$ and $2g+n\ge 12$ (so in particular it is false for $\Mct_6$).

So the Gorenstein conjecture for smooth curves discussed in this paper was the last of the three to fall, and the counterexamples to the other Gorenstein conjectures gave reason to be skeptical of this one.

\subsection{Faber-Zagier and $3$-spin relations}

Faber's original method \cite{Faber} for constructing tautological relations on $\M_g$ did not produce a family of relations that was easy to describe. However, Faber and Zagier looked at the relations produced by this method in small genus and guessed a relatively simple formula (involving hypergeometric generating functions) giving a finite list of kappa polynomials for each $g,d$ that seemed to span the same space of relations. These are now known as the {\it Faber-Zagier relations}.

The Faber-Zagier relations were first proved to vanish (i.e. to be actual relations) via moduli of stable quotients by Pandharipande and the author in \cite{PP-FZ}, in which we also proposed the following conjecture:
\begin{fzconj}[{\cite[Conjecture 2]{PP-FZ}}]
  The Faber-Zagier relations span the space of tautological relations on $\M_g$.
\end{fzconj}
This conjecture agrees with the Gorenstein conjecture for $g\le 23$ (and thus is true there) but contradicts it for (seemingly all) $g \ge 24$. For $g=24,25,26$ exactly one of the two conjectures is true (since their predictions differ by a single relation), but for $g\ge 27$ it is possible that neither conjecture is true.

So Theorem~\ref{thm:main} verifies the Faber-Zagier conjecture for $g=25$, but it is still open for $g=24$ and all $g\ge 26$. All known methods to prove and compute tautological relations on $\M_g$ for $g\ge 24$ have only produced relations in the span of the Faber-Zagier relations. This failure to find relations beyond the Faber-Zagier relations was weak evidence in favor of this conjecture and against the Gorenstein conjecture.

The $3$-spin relations \cite{PPZ,Janda} are a generalization/extension of the Faber-Zagier relations from $\M_g$ to $\Mbar_{g,n}$; they are a formula giving a large finite list of (proven) relations in every $\R^d(\Mbar_{g,n})$. Taking $n=0$ and restricting to $\M_g$ recovers the Faber-Zagier relations exactly. The author has conjectured \cite{Pixton-conj} that the $3$-spin relations (restricted to open loci if necessary) span all tautological relations in $\Mbar_{g,n},\Mct_{g,n}$, and $\Mrt_{g,n}$. For the purposes of this paper, we will only need the restriction of the $3$-spin relations to $\M_{g,1}$, i.e. an analogue of the Faber-Zagier relations with one marked point. (We could alternatively use the tautological relations constructed by Yin \cite{Yin} using the universal Jacobian over $\M_{g,1}$, since they have been computed to give the same space of relations for the relevant genera.)

\section{Searching for a detector}\label{sec:computations}
Let us explain how to carry out the plan in the introduction for proving Theorem~\ref{thm:main}. We will describe the approach for general $g,d$; our goal is to detect that a class $P\in\R^d(\M_g)$ is nonzero. (So for Theorem~\ref{thm:main}, we have $g=25,d=12$.) We want to construct a detector class
\[
  Q = \lambda_g\lambda_{g-2}\alpha\in \R^{3g-3-d}(\Mbar_g),
\]
where $\alpha\in\R^{g-1-d}(\Mbar_g)$. We want $Q$ to vanish on the entire boundary of $\Mbar_g$. Since $\lambda_g$ already vanishes on the divisor of stable curves with a non-separating node, we just need to consider pullbacks $j_h^*Q$ for gluing maps
\[
  j_h:\Mbar_{h,1}\times\Mbar_{g-h,1}\to \Mbar_g,\quad\quad 1\le h\le\lfloor\frac{g}{2}\rfloor.
\]
In addition, the $\lambda_g\lambda_{g-2}$ factor constrains which tautological classes $\alpha$ are potentially useful to us. The only nontrivial dual graphs that survive have a single separating edge. So we will only consider $\alpha$ of the form
\[
  \alpha = A + \sum_{h=1}^{\lfloor\frac{g}{2}\rfloor}(j_h)_*A_h,
\]
where $A(\kappa_i)$ is a kappa polynomial of degree $g-1-d$ and $A_h(\kappa_i^{(1)},\psi_1,\kappa_i^{(2)},\psi_2)$ is a polynomial of degree $g-2-d$ in the kappa and psi classes pulled back from the factors of $\Mbar_{h,1}\times\Mbar_{g-h,1}$.

When we pull back along $j_h$, the $\lambda_g\lambda_{g-2}$ factor becomes $\lambda_h\lambda_{h-1}\boxtimes\lambda_{g-h}\lambda_{g-h-1}$ and thus vanishes on the boundary of $\Mbar_{h,1}\times\Mbar_{g-h,1}$. As a consequence, it suffices to find $\alpha$ such that $(j_h^\circ)^*\alpha = 0$, where $j_h^\circ:\M_{h,1}\times\M_{g-h,1}\to\Mbar_g$ is the restriction to the interior of $j_h$.

Let $\kappa_i^{(1)},\psi_1,\kappa_i^{(2)},\psi_2$ be the kappa and psi classes on the factors of $\M_{h,1}\times\M_{g-h,1}$. Then we have
\[
  (j_h^\circ)^*\alpha = A|_{\kappa_i := \kappa_i^{(1)}+\kappa_i^{(2)}} + (-\psi_1-\psi_2)A_h.
\]

Given lists of tautological relations on $\M_{h,1},\M_{g-h,1}$ such as the restricted $3$-spin relations (even if we don't know that they are the full lists of relations), we can now search for polynomials $A,A_h$ such that
\begin{equation}\label{eq:BVcond}
  A|_{\kappa_i := \kappa_i^{(1)}+\kappa_i^{(2)}} + (-\psi_1-\psi_2)A_h = 0
\end{equation}
follows from those relations (for each $1 \le h \le \lfloor\frac{g}{2}\rfloor$). (This needs to be modified slightly if $g$ is even and $h=\frac{g}{2}$, but we are primarily interested in the case $g=25$ so we will ignore this for simplicity.)

If we successfully find such $A,A_h$, then we know the resulting $Q$ has the desired boundary vanishing property. We aren't done yet - we also need to compute the pairing of $Q$ against $P$. For the $A_h$ terms, this only requires computing the $\lambda_g\lambda_{g-1}$ Hodge integrals that were already needed for the proportionalities in the Gorenstein conjecture. But for the $A$ term, we need to compute Hodge integrals of the form
\[
  \int_{\Mbar_g}\lambda_g\lambda_{g-2}K,
\]
where $K$ is a kappa polynomial. This is harder but is computationally feasible in a couple different ways, e.g. via a recursion described by Getzler-Pandharipande \cite[Theorem 4]{Getzler-Pandharipande}.

\subsection{Implementation}
The existence of a $Q_{25}$ with the desired properties (boundary vanishing and nonzero pairing with $P_{25}$) was originally checked by code written by GPT-6 Astra. Instead of reading Astra's code and trying to verify that it works as intended, the author wrote code from scratch (making use of the admcycles Sage package \cite{admcycles} for some things) to perform the same search. This code can be found at \url{https://websites.umich.edu/~pixton/M25/M25.sage}. We briefly describe here how it works.

The admcycles package is used to compute all the $3$-spin relations on $\M_{h,1}$ for $1\le h\le g-1$. Quotienting the polynomial algebra generated by formal kappa and psi variables by these relations gives graded algebras $\mathcal{Q}^*(\M_{h,1})$.

The data $(A,(A_h))$ is then interpreted as an element of the $\QQ$-vector space
\[
  \QQ[\kappa_i]_{\text{deg $(g-1-d)$}} \oplus \left(\bigoplus_{h=1}^{\lfloor\frac{g}{2}\rfloor}\mathcal{Q}^*(\M_{h,1})\otimes_\QQ\mathcal{Q}^*(\M_{g-h,1})\right)_{\text{deg $(g-2-d)$}}.
  \]
  The boundary vanishing conditions \eqref{eq:BVcond} can then be interpreted as defining a linear map to the similar $\QQ$-vector space
  \[
    \left(\bigoplus_{h=1}^{\lfloor\frac{g}{2}\rfloor}\mathcal{Q}^*(\M_{h,1})\otimes_\QQ\mathcal{Q}^*(\M_{g-h,1})\right)_{\text{deg $(g-1-d)$}},
  \]
  and elements of the kernel of this map yield classes
  \[
    Q := \lambda_g\lambda_{g-2}\left(A + \sum_{h=1}^{\lfloor\frac{g}{2}\rfloor}(j_h)_*A_h\right)
    \]
    with the desired boundary vanishing property.

In the case of $g=25$, these two vector spaces have dimensions $8805$ and $8738$, and the map turns out to be surjective, so the kernel has dimension $67$. This includes examples $(A,(A_h))$ given by restricting a $3$-spin relation on $\Mbar_{25}$ to the appropriate strata. After removing those trivial examples (which form a $33$-dimensional family), there is a $34$-dimensional space of candidates for $Q_{25}$.

It remains to compute the pairing of these $Q_{25}$ candidates with arbitrary kappa polynomials. The admcycles package is used to compute the $\lambda_g\lambda_{g-1}$ integrals that show up. The $\lambda_g\lambda_{g-2}$ integrals are harder. For ease of implementation, the author didn't try to implement the recursion \cite[Theorem 4]{Getzler-Pandharipande} (though GPT-6 Astra claimed to have used this approach successfully). Instead, the author precomputed a table of all integrals of pure kappa monomials on $\Mbar_g$ and combined that with an explicit identity expressing $\lambda_g\lambda_{g-2}$ as a kappa polynomial \cite[Theorem 1.2]{ABDKS}. Although perhaps more computationally demanding, this was more straightforward to implement and worked fine for $g\le 25$.

Combining this pairing with the regular socle pairing, the resulting rank is $44$, one higher than the rank of the regular socle pairing alone. Therefore the dimension of $\R^{12}(\M_{25})$ is at least one higher than predicted by the Gorenstein conjecture, establishing Theorem~\ref{thm:main}. (In fact $\dim_{\QQ}\R^{12}(\M_{25}) = 44$ because the Faber-Zagier relations give that as an upper bound.)

\subsection{Additional computations}
Of course, there is nothing preventing one from trying the same approach to prove nonvanishing of other ``missing relations''. The linked code can also attempt this approach for the missing relation in genus $24$ (also in codimension 12), but it does not succeed. In that case, the only $(A,(A_h))$ satisfying the linear conditions (as implemented there) are the restrictions of $3$-spin relations.

GPT-6 Astra's better optimized code was able to attempt this same approach through genus $30$ (and could certainly go further with a little more computing time). For $g = 26,28,30$ and $d=\frac{g}{2}$ things are the same as for $g=24$: there are no nonvanishing $Q$ candidates. On the other hand, for $g=27,29$ and $d=\frac{g-1}{2}$ things are much the same as for $g=25$: the approach is successful at proving that the ring is not Gorenstein and that the Faber-Zagier relations span all relations in $\R^d(\M_g)$. However, this does not fully determine the tautological ring in these cases because there are missing relations in other codimensions too - for instance the method fails to prove nonvanishing of a missing relation in $\R^{14}(\M_{27})$.

It seems like this odd genus vs even genus discrepancy is really just a consequence of the approach being less effective for proving nonvanishing of classes in $\R^d(\M_g)$ for high $d$. More precisely, it seems like nonvanishing $Q$ (divisible by $\lambda_g\lambda_{g-2}$) with the boundary vanishing property only exist when $d\le\frac{g-1}{2}$, but missing relations only seem to exist when $d\ge\frac{g-1}{2}$. So there is (experimentally) only a single location $d=\frac{g-1}{2}$ where this approach is effective.

Finally, the author asked GPT-6 Astra to write code to implement this method for rational tails spaces $\Mrt_{g,n}$ in various minimal cases where the Faber-Zagier conjecture and the Gorenstein conjecture are known to be contradictory, as described in \cite[Conjecture A.1]{Pixton-thesis} (with rational tails instead of powers of the universal curve). Astra claims that the method works in those cases (e.g. in $\R^{10}(\M_{20,1})$ and in $\R^9(\Mrt_{17,2})$), but the author has not yet verified or independently implemented the code.

\section{Further questions}\label{sec:questions}
Some further questions are obvious - we still do not know for precisely which values of $g$ the Gorenstein conjecture is true, or whether the Faber-Zagier conjecture might be true for all $g$ (or for any $g$ besides $g\le 23$ and $g=25$). 

In this section we only list a few questions relating to the specific $\lambda_g\lambda_{g-2}$ approach used in this paper. First, a question about the space of tautological classes that can be used for this approach:
  
\begin{question}
  Let $V_{g,d}\subseteq \QQ[\kappa_i]_{\text{deg $(g-1-d)$}}$ be the subspace of kappa polynomials that can be used as $A$ in this approach (i.e. admit tautological boundary corrections such that boundary vanishing can be certified by $3$-spin relations on $\M_{h,1}$). We know that $V_{g,d}$ contains all Faber-Zagier relations, but it is often strictly larger. Is there some other characterization of $V_{g,d}$?
\end{question}

And a series of questions specifically about odd genus, motivated by numerical experimentation up through genus $25$:
\begin{question}
  Let $m\ge 1$.
  \begin{enumerate}[(a)]
  \item Is it always true that $\kappa_m\in V_{2m+1,m}$? This is equivalent to the statement that $\kappa_m\in \R^*(\M_{2m+1-i,1})$ is divisible by $\psi^i$ for all $1\le i \le m$.
  \item Assuming (a), is it always true that there is a unique choice of kappa/psi polynomials $A_h$ (modulo relations) that work with $A = \kappa_m$? If so, can we interpret the resulting cycle in some way or give a formula for the $A_h$?
  \item Assuming (a) and (b), is it true for all $m\ge 12$ that the corresponding $\lambda_g\lambda_{g-2}$ functional on $\R^m(\M_{2m+1})$ is not in the span of the $\lambda_g\lambda_{g-1}$ functionals (and thus the Gorenstein conjecture is false for $g=2m+1$)?
    \end{enumerate}
  \end{question}

\bibliographystyle{plain}
\bibliography{M25}

\vspace{8 pt}
\noindent Department of Mathematics, University of Michigan\\
pixton@umich.edu

\end{document}